\documentstyle[12pt]{article}
\begin{document}

\vspace{1.cm}

\begin{center}
{\huge {\bf Paul levy's Isoperimetric Problems on cyclic polygons}} \footnote{
Mathematics Subject Classification 51M10, 51M25, 52A40} \\

\end{center}

\vspace{2.cm} {\large \ Abd Raouf Chouikha} \footnote{University of Paris-Nord\
LAGA,CNRS UMR 7539,\\ 
{\small chouikha@math.univ-paris13.fr}}

\vspace{3.cm}

{\bf Abstract}\newline
{\it In this paper we are interested in isoperimetric
inequalities for plane n-gons in relation with a old conjecture proposed
by P. Levy.}\newline

\vspace{2.cm}

\section{Introduction}

The problem of isoperimetry is often associated with the work of Zenodorus and Pappus in the ancient world.
This old problem in geometry has a significant impact on various branches of mathematics. The problem deals with finding a
closed curve or surface that encloses the maximum area or volume for a given perimeter or boundary length. It was not until the work of Steiner (1838), who provided a rigorous proof, \cite{b}.\\

This problem can be expressed as follows:, for a simple closed curve ${\cal C}$ (in the
euclidian plane) of length $L$ enclosing a domain of area A, the following
inequality holds

\begin{equation}  \label{1}
L^2 - 4\pi A \geq 0.
\end{equation}
Equality is attained if and only if this curve is a euclidean circle. This
means that among the set of domains of fixed area, the euclidean circle has
the smallest perimeter.\newline
 Steiner's proof was completed later by several other mathematicians. Notably by Paul Levy, who introduced and popularized symmetrization processes (akin to Steiner or Minkowski symmetrization) applied to curves and polygons.\\
In 1902, Hurwitz published a short proof using the Fourier series that applies to piecewise closed curves

The above inequality (\ref{1}) could be easily deduced from the Wirtinger
inequality

\begin{equation}
\int^{2\pi}_0|f^{\prime}(x)|^2 dx \geq \int^{2\pi}_0|f(x)|^2 dx,  \label{2}
\end{equation}
where $f(x)$ is a continuous periodic function of period $2\pi$ whose
derivative $f^{\prime}(x)$ is also continuous and $\int^{2\pi}_0 f(x) dx = 0$%
. Equality holds if and only if $f(x) = \alpha \cos x + \beta \sin x$ (see
Osserman \cite{o1} or Berger \cite{b}).\newline
For any curve ${\cal C}$ of length $L$ enclosing an area $A$, the quantity
\quad $L^2 - 4\pi A$\quad is called {\it the isoperimetric deficiency} of\ $%
{\cal C}$,\ because it decreases towards zero when \ ${\cal C}$ \ tends to a
circle. \newline

 It is well-known that these inequalities are also useful in
mathematical physics statistics and other applied sciences.
More precisely, it is related to the principle of least action in physics, in that it can be expressed: what is the principle of action which encloses the greatest area, with the greatest economy of effort.\\

\subsection{Connexion with the Bonnesen inqualities }

A Bonnesen inequality introduces an additional positive "defect" term on the right side, usually depending on the radii of the inscribed circle ($r$) and circumscribed circle ($R$). It quantifies exactly how much a shape deviates from being a perfect circle. When applied specifically to $n$-gons (polygons with $n$ sides), the boundary is constrained, and the continuous circle is replaced by a regular $n$-gon as the optimal, area-maximizing shape.

As an extension, Bonnesen proves \cite{o1} that if  ${\cal C}$ is convex and
there exists a circular annulus containing ${\cal C}$ of thickness $d$, then we
have
$$4 \pi d^2 \leq L^2 - 4\pi A.$$

Let $K$ be the region bounded by ${\cal C}$. The radius of the smallest circular disk containing $K$ is called the
circumradius, denoted $R$  . The radius of the largest circular disk contained in $K$ is the inradius, denoted $r$.
There is a related isoperimetric inequality known as the Bonnesen :

\begin{equation}
L^2 - 4\pi A \geq \pi^2 (R - r)^2,  \label{3}
\end{equation}where\ $R$ \ is the circumradius and \ $r$ \ is the inradius of the curve \ $%
{\cal C}$. \newline
Note that if the right side of (\ref{3})  equals zero, then $R = r$.
This means that ${\cal C}$ is a circle and $L^2 - 4\pi A = 0$.\newline
More generally, inequalities of the form 
\begin{equation}
L^2 - 4\pi A \geq K,  \label{4}
\end{equation}
are called Bonnesen-type isoperimetric inequalities if  equality is only
attained for the euclidean circle. In the other words, $K$ is positive and
satisfies the condition

\[
K=0 \quad \mbox{ implies } \quad \ L^2 - 4\pi A = 0.
\]
(See \cite{o2} for a general discussion and different generalisations).\newline

\subsection{Inequalities for $n$-gons}
Let a $n$-gon  with side of lenght $a_i$, its perimeter is \ $L_n  = a_1 + a_2 +...+a_n.$\ It is known a necessary condition for the sides 
$a_i, \ 1<i<n$  to be the edge lengths of a cyclic polygon is
that each of the $a_i$ is less than the sum of the rest of them. Or equivalently
\begin{equation} a_n < a_1 +a_2 +...+ a_{n-1}, \Leftrightarrow \frac{2 a_i}{L_n} < 1.\end{equation} 
 It is alo proved that then there exists a unique (up to an isometry) convex cyclic polygon with sides lengths $a_1 ,...,a_n$,\ \cite{b} or \cite{p}, Theorem 1. Moreover, in studying geometric inequalities for polygons, one only
needs to consider cyclic polygons. Indeed, such polygons always enclose the largest area.\\

We also need the following, \cite{o2}\\

{\bf Lemma 1.}\ {\it Let $K$ be a domain bounded by a rectifiable Jordan curve ${\cal C}$ with $L, A, R, r$\ its length,
area, circumradius and inradius. Then there exists a sequence of polygons $\Pi_n$ that
$$A n \rightarrow A,\quad L n \rightarrow L,\quad and\quad R_n   \rightarrow R  .$$
If in addition $K$ is convex, we may also arrange that\quad  $r_n \rightarrow r$  .}\\

 That means, intuitively, we may approximate ${\cal C}$ by finer polygonal
curves.\\

For any $n$-gon $\Pi_n$ with perimeter $L_n$ and area $A_n$, the baseline inequality (analogous to the circular isoperimetric inequality) is:

$$L_n^2 - 4(n \tan{\frac{\pi}{n}}) A_n \geq 0,$$ 

Equality holds if and only if the $\Pi_n$ is regular.\\

To turn this into a Bonnesen-style inequality, we introduce geometric radii. For a polygon, we look at its inradius $r$ (the largest circle contained inside) and circumradius $R$ (the smallest circle containing it).

For general or convex $n$-gons, several classic bounds exist that introduce these radii to bound the isoperimetric deficit:

A direct relationship involving the inradius $r$ is:
$$A_n \leq L_n r - nr^2 \tan{\frac{\pi}{n}}.$$
Rearranging this quadratic inequality in terms of its discriminant is precisely what yields the strong polygonal Bonnesen inequality, \cite{z}:\begin{equation} L_n^2 - 4n \tan\left(\frac{\pi}{n}\right) A_n \geq \left[ L_n - 2n r \tan\left(\frac{\pi}{n}\right) \right]^2 = \left[L_n - 2n R \sin\left(\frac{\pi}{n}\right)\right]^2.\end{equation}

Notice that (6) was improved by \cite{zd}, Corollary 6 \begin{equation}L_n^2 - 4n \tan \frac{\pi}{n}\ A_n \geq 2 R \tan \frac{\pi}{n}\left[ L_n - 2n R \sin \frac{\pi}{n} \right].\end{equation}

If we look at the difference between the circumradius $R$ and inradius $r$ (this one often called annulus width), a classic Bonnesen inequality for a cyclic polygon can always also be expressed using circular radii.
 $$L^2 - 4 \pi A \geq \pi^2 (R-r)^2.$$

On the other hand, in the Hausdorff metric space of convex bodies, a standard Bonnesen inequality provides a quantitative stability estimate. For a general convex body $K$, it tells us how close $K$ is to a ball.For the space of polygons with $n$ sides ($\mathcal{P}_n$), the "optimal" shape (the one that minimizes the isoperimetric quotient) is the regular $n$-gon, $P_n^*$. The polygonal Bonnesen inequality provides a quantitative measure of how close a given polygon $P \in \mathcal{P}_n$ is to $P_n^*$.If we define the isoperimetric deficit of an $n$-gon as:
$$\Delta_n(P) = L^2 - 4n \tan\left(\frac{\pi}{n}\right) A$$Then a sharp Bonnesen-style stability theorem states:$$\Delta_n(P) \geq c_n \cdot \delta_H(P, P_n^*)^2$$where $\delta_H$ is the Hausdorff distance (up to scaling and translation) and $c_n$ is a positive constant depending purely on $n$. This proves that if the deficit is small, the polygon is geometrically close to being regular.\\

Moreover, we know that an $n$-gon has a maximum area among all $n$%
-gons with the given set of sides if it is convex and inscribed in a circle.
Let \quad $a_1, a_2, ....,a_n$ \quad denote the lengths of the sides of \ $
\Pi_n.$ \ For a triangle, Heron's formula gives the area

\[
A_3 = {\frac{1}{4}} L_3^2 \sqrt{ (1 - {\frac{2a_1}{{L_3}}})(1 - {\frac{2a_2}{%
{L_3}}})(1 -{\frac{2a_3}{{L_3}}})}.
\]
For a quadrilateral, the Brahmagupta formula gives a bound for the area

\[
A_4 \leq {\frac{1}{4}} L_4^2 \sqrt{ (1 - {\frac{2a_1}{{L_4}}})(1 - {\frac{%
2a_2}{{L_4}}})(1 - {\frac{2a_3}{{L_4}}})(1 - {\frac{2a_4}{{L_4}}})}
\]
with equality if and only if \ $\Pi_4$ \ can be inscribed in a circle.

\subsection{ Isoperimetric constants}

We can ask if it is possible to get an analogous formula for other plane  polygons (not 
necessarily inscribed in a circle). More precisely, is
the area $A_n$ of the $n$-gon is close to the following expression ?

\begin{equation}
P_n = {\frac{ L_n^2}{{4}}} \sqrt{ (1 - {\frac{2a_1}{{L_n}}})(1 - {\frac{2a_2%
}{{L_n}}})(1 - {\frac{2a_3}{{L_n}}}).....(1 - {\frac{2a_n}{{L_n}}})}
\label{6}
\end{equation}
called the pseudo-area.\\

This question has been considered by many geometers who tried to
compare\quad $\displaystyle A_n$\quad with \quad $P_n$. Among them Robbins \cite{rob}, Zhang \cite{z}, \cite{kkz} , Dulio-Laeng \cite{dl}, Pak \cite{pak}, Svrtan \cite{sv}, Pech \cite{pech},... The one of the first among them is undoubtedly Paul Levy, who addressed this problem in his 1966 article \cite{l} and even proposed a conjecture that remains relevant today. More precisely, he defined 
 the ratio $${\displaystyle \varphi_n = {\frac{A_n}{{P_n}}}},$$
for any $n$-gon \ $\Pi_n$, with sides $a_1,a_2,....,a_n$
enclosing an area $A_n,$ and $P_n$ defined as above. However, he remarked for the regular $n$-gon

\begin{equation}
{\varphi_n}^0 = \frac{{A_n}^0}{{P_n}^0} = {\frac{1}{{{n \tan{\frac{\pi}{{n}}}}\ (1-{\frac{2}{{n}}}%
)^{n/2}}}} \label{7}
\end{equation}
He noticed that \ ${\varphi_n}^0$\ is decreasing as a function of $n$ and should verify   
$${\varphi_3}^0={\varphi_4}^0=1,\qquad \frac{e}{\pi} < {\varphi_n}^0 < 1.$$
Moreover, Paul Levy proposed the following\\ 

{\bf Conjecture (L):} {\it Define the ratio ${\displaystyle \varphi_n = {\frac{%
A_n}{{P_n}}}}$. For any $n$-gon \ $\Pi_n$, with sides $a_1,a_2,....,a_n$
enclosing an area $A_n,$ and $P_n$ defined as above, this ratio verifies }
$${\varphi_n}^0 < \varphi_n < 1.$$\\

In this paper, we revisit this conjecture and examine the underlying problems using concrete cases that are more general than those presented by Paul Levy, thereby allowing for a better understanding of the issue at hand. We hope that the method we employ will lead to a complete resolution of the conjecture. More precisely, we consider an $(n+1)$-gon in which $n-1$ sides are equal to $a$, while the remaining two have lengths $b$ and $c$. We also consider a $(k+m)$-gon with $k$ sides equal to $a$ and $m$ sides equal to $b$, thus generalizing the Macnab polygon case (where $k=m$) used in \cite{ch2}. We also improve Proposition 4-2 of \cite{ch1} in considering the case of $n$-gon whose sides all satisfy the condition \ $a_i < \frac{k\ L_n}{n}, \ 1<i<n,\ k >1$.\ Further examples and additional approaches will be explored with the aim of resolving this thorny conjecture. \newline

\section{isoperimetric inequalities for some polygons}

\subsection {Curved $n$-gons}
In this section, we will re-examine and develop examples provided by Paul Levy \cite{l} to illustrate his conjecture. Notice that he has only considered the case \ $n \rightarrow \infty$.

\subsubsection{First example}
Let $AB$ be an arc of a circle of length $2\alpha$ and radius $R=1$. consider its contour as the boundary of a polygon with $(n+1)$-sides, one of which is the chord $A$B of length $a = 2sin \alpha$, and the others subtending arcs of the same length $2a/n$. Let $L_{n+1}, A_{n+1},$ and $P_{n+1}$ be respectively its perimeter, its area and pseudo-area (5) :
$$L_{n+1} = 2\,n\sin  {\frac {\alpha}{n}}  +2\,\sin  \alpha
 , \qquad A_{n+1} = \frac{n}{2}\sin {\frac {2\alpha}{n}}  -\frac{1}{2} \sin  2\,\alpha ,$$
$$P_{n+1} = \frac{L_{n+1}^2}{4} \left( 1-\frac{4\,\sin  {\frac {\alpha}{n}} }{ {L_{n+1}}} \right) ^{\frac{n}{2}}\sqrt {1-4\,{\frac {\sin  \alpha  }{
L_{n+1}}}} =$$ $$ \left( \alpha+\sin  \alpha   \right) ^{2} \left( 1-\frac{2\,
\sin \left( {\frac {\alpha}{n}} \right)}{  \left( \alpha+\sin  
\alpha   \right)} \right) ^{\frac{n}{2}}\sqrt {1-2\,{\frac {
\sin \left( \alpha \right) }{\alpha+\sin  \alpha  }}}.$$
The quotient $\phi_{n+1} = \frac{A_{n+1}}{P_{n+1}}$ is
$$\phi_{n+1} = \frac{ \left( \frac{n}{2}\sin \left( 2\,{\frac {\alpha}{n}} \right) -\frac{1}{2}
\sin  2\,\alpha \right)  }{  \left(  \left( 
n\sin \left( {\frac {\alpha}{n}} \right) +\sin  \alpha  
 \right) ^{1-\frac{n}{2}} \right) \left(  \left( n\sin \left( {\frac 
{\alpha}{n}} \right) +\sin  \alpha  -2\,\sin \left( {
\frac {\alpha}{n}} \right)  \right) ^{\frac{n}{2}} \right)}{\frac {1}{
\sqrt {{n}^{2} \left( \sin  {\frac {\alpha}{n}}  
 \right) ^{2}- \left( \sin  \alpha   \right) ^{2}}}}.$$
Following \cite{l} as limit when $n \rightarrow \infty$ we find  
$$L(\alpha) = 2 \alpha + 2 \sin \alpha, \ A(\alpha) = \alpha - \sin \alpha \cos \alpha,$$
$$\phi(\alpha) = \frac{ \left( \alpha-\frac{1}{2}\sin  2\,\alpha  
 \right) {{\rm e}^{{\frac {\alpha}{\alpha+\sin  \alpha  }
}}}}{ \left( \alpha+\sin  \alpha   \right)}{
\frac {1}{\sqrt {{\alpha}^{2}- \left( \sin  \alpha  
 \right) ^{2}}}}.$$
One proved \cite{ch1}, Proposition 4-1, for\  $0 <\alpha <\pi$ \ the function\  $\phi(\alpha) $\ is strictly decreasing and\  $\sqrt{\frac{e}{3}} > \phi(\alpha) > \frac{e}{\pi}.$\\
For  \ the $(n+1)$-gon,\ and for \ $n =2, 3$\ one finds $$\phi_3 = -{\frac {\sqrt {4} \left( -1+\cos  \alpha   \right) 
}{4\sqrt { \left( 1- \left( \cos  \frac{\alpha}{2}   \right) ^{
2} \right) ^{2}}}} =1,$$ $$ \phi_4 = {\frac { 2\left( \sin \frac{2\alpha}{3}  \right) ^{3}\sqrt 
{3\,\sin  \frac{2\alpha}{3}  +\sin  \alpha  }}{
 \left( \sin \frac{\alpha}{3} +\sin  \alpha  
 \right) ^{\frac{3}{2}}\sqrt {9\, \left( \sin \frac{\alpha}{3} 
 \right) ^{2}- \left( \sin  \alpha   \right) ^{2}}}} =1.$$
For \ $n > 3$\ and \ $0< \alpha < \pi$\ we may verify by {\it Maple} that \ $\phi_{n+1}$ \ is strictly decreasing with respect to $n$ and 
$$\phi_{n+1}^0  < \phi_{n+1} <1.$$ Conjecture (L) is then verified.\\ 

\subsubsection{Second example} 
Consider on the same circle of center $O$ and radius $1$ a variable point $C$ on the arc $AB$ opposite the one enclosing area $A(\alpha)$. Let $2\theta$ denote the angle $AOC$,\ $A_{n+2}(\alpha, \theta)$ the sum of $A_{n}$ and the area of ​​triangle $ABC$. We find the length, the area and the pseudo-area of this polygon with $(n+2)$ sides:
$$L_{n+2} = 2\,n\sin \left( {\frac {\alpha}{n}} \right) +2\,\sin \left( \alpha+
\theta \right) +2\,\sin  \theta , \ A_{n+2} = \frac{n}{2}\sin \left( {\frac {2\alpha}{n}} \right) -\sin  \alpha
  \cos \left( \alpha+2\,\theta \right),$$
$$P_{n+2} = \frac{{L}^{2}}{4} \left( 1-\frac{2\,\sin \left( {\frac {\alpha}{n}} \right)}{ {L}
} \right) ^{\frac{n}{2}}\sqrt {1-{\frac {2\sin  \alpha  
}{L}}}\sqrt {1-{\frac {2\sin \left( \alpha+\theta \right) }{L}}}.$$ The quotient is
$$\phi_{n+2} = \frac{4\, \left( \frac{n}{2}\sin  {\frac {2\alpha}{n}}  -\sin
  \alpha  \cos \left( \alpha+2\,\theta \right)  \right)}{ {
L}^{2}   \left( 1-\frac{\sin  {2\frac {\alpha}{n}} }{ {L
}} \right) ^{\frac{n}{2}}  \sqrt {1-{\frac {2
\sin  \alpha  }{L}}}{\sqrt {1-2\,{\frac {\sin
 \left( \alpha+\theta \right) }{L}}}}}.$$

For infinite $n$, we have shown that $\phi(\alpha,\theta)$ has a maximum for \ $\theta =\frac{\pi-\alpha}{2},$ \ \cite{ch1}. It can easily be considered that for sufficiently large \ $n,\ \phi_{n+2}$\  admits a maximum for \ $\theta =\frac{\pi-\alpha}{2}.$ In this case a calculation yields 
$$\phi_{n+2}(\alpha,\theta) < \phi_{n+2}(\alpha,\frac{\pi-\alpha}{2}) =$$ $$ \frac{ \left( \frac{n}{2}\sin \left( {\frac {2\alpha}{n}} \right) +\sin \alpha   \right)}{  \left( n\sin  {\frac {
\alpha}{n}} + 2\,\cos  \frac{\alpha}{2}  -\sin  {
\frac {\alpha}{n}}   \right) ^{\frac{n}{2}}  }\frac {\left( n\sin  {\frac {\alpha}{n}}  +2\,\cos
  \frac{\alpha}{2}   \right) ^{\frac{n}{2}-1}}{
\sqrt {n\sin  {\frac {\alpha}{n}}  +2\,\cos  \frac{
\alpha}{2}  -\sin  \alpha  }\ {\sqrt {n\sin
  {\frac {\alpha}{n}}  +\cos  \frac{\alpha}{2}  }}}
     .$$
		
	The limit case ( $n \rightarrow \infty$ ) is 
	$$ \frac{\left( \alpha+\sin  \alpha   \right) {{\rm e}^{{\frac 
{\alpha}{\alpha+2\,\cos  \frac{\alpha}{2}  }}}}}{\alpha \left( \alpha+
2\,\cos \frac{\alpha}{2} \right) },$$
which is strictly decreasing from $1$ to $\frac{e}{\pi}$.

\subsection{Special polygons} \ {\bf 1.}\ Let us consider a regular $n$-gon with sides of lenght $a$. We will replace one side by two sides of same lenght $b$. We then obtain a $(n+1)$-gon whose $n-1$ sides are of lenght $a$ and two sides of lenght $b$. An easy computation yields its perimeter and its area $$b =  a {\sin \frac{\pi}{2n}}{\sin \frac{\pi}{n}}  = \frac{a}{2 \cos \frac{\pi}{2n}},\qquad
L_{n+1} = (n-1) a + \frac{a}{ \cos \frac{\pi}{2n}} ,$$
$$A_{n+1} = \frac{n a^2}{4} \cot \frac{\pi}{n} + \frac{ a^2}{4} \left(\frac{1}{\sin \frac{\pi}{n}} - \cot \frac{\pi}{n}\right) =\frac{n a^2}{4} \cot \frac{\pi}{n} + \frac{ a^2}{4}  \frac{1}{\sin \frac{\pi}{n}}  \left(1 - {\cos \frac{\pi}{n}}\right)=$$
$$\frac{1}{4} \frac{a^2 \left((n-1) \cos \frac{\pi}{n} +1 \right)}{\sin \frac{\pi}{n}}.$$

The function $P_{n+1}$ given by (5) may be computed
$$P_{n+1} =  {\frac{ L_{n+1}^2}{{4}}} \prod_{1<i<n+1}{\left(1 - {\frac{2a_i}{{L_n}}}\right)^{\frac{1}{2}}} = {\frac{ L_{n+1}^2}{{4}}}\sqrt { \left( 1-{\frac {2 a}{L_{{n+1}}}}
 \right) ^{n-1} \left( 1-{\frac {2 b}{L_{{n+1}}}} \right) ^{2}} =$$ $$ {\frac{ L_{n+1}^2}{{4}}} \left( 1-{\frac {2 a}{L_{{n+1}}}} \right) ^{\frac{n-1}{2}} \left( 1-{\frac{2\,a\sin {\frac {\pi }{2n}}}{ \sin  {\frac {\pi }{n}} \   {L_{{n+1}}}}} \right).$$
The quotient $\phi_{n+1} = \frac{A_{n+1}}{P_{n+1}} $ can also be computed
$$\phi_{n+1} =  \frac{a^2 \left((n-1) \cos \frac{\pi}{n} +1 \right)}{\sin \frac{\pi}{n} \ L_{n+1}^2  \left( 1-{\frac {2 a}{L_{{n+1}}}} \right) ^{\frac{n-1}{2}} \left( 1-\frac{2\,a\sin {\frac {\pi }{2n}}}{  
  \sin  {\frac {\pi }{n}} \  {L_{{n+1}}
}} \right)}=$$

$$\frac{ \left((n-1) \cos \frac{\pi}{n} +1 \right)}{\sin  {\frac {\pi }{n}}   \left(  \left( n-1 \right) +\frac{1}{
  \cos  {\frac {\pi }{2n}} } 
 \right) ^{2} \left( 1-\frac{2 }{   \left( n-1 \right) +\frac{1} {\left( \cos
 {\frac {\pi }{2n}} \right)    }}
 \right) ^\frac{n-1}{2} \left( 1-\frac{2\,\sin {\frac {\pi }{2n}}}{
   \left( \sin  {\frac {\pi }{n}}   \right) 
 \left(  \left( n-1 \right) +\frac{1}{ \left( \cos {\frac {\pi }{
2n}}   \right)}  \right) } \right) }.$$

As we can see, $\phi_{n+1}$ is independent of the length $a$. For $n=3$ and $n=2$ one finds of course $$\phi_{4} = \frac{\frac{4}{3}\,\sqrt {3}  \left( 1-
2/3\,{\frac {\sqrt {3}}{2+\frac{2}{3}\sqrt {3}}} \right) }{ \left( 2+\frac{2}{3}\sqrt {3} \right) ^{2}
 \left( 1-{\frac {2}{2\,+\frac{2}{3}\sqrt {3}}} \right)}  =1,$$
$$\phi_{3}= \frac{1}{ \left( 1+\sqrt {2} \right) ^{2}}{\frac {1}{\sqrt {1-{
\frac {2}{1+\sqrt {2}}}}}} \left( 1-{\frac {\sqrt {2}}{1+\sqrt {2}}
} \right)  =1.$$

By {\it Maple} it is verified that $\phi_{n+1}$ is strictly decreasing with respect to $n$ from $1$ to $$lim_{n\rightarrow\infty} \phi_{n+1} = \frac{e}{\pi}.$$

Moreover, comparing $\phi_{n+1}$ with the analog for the regular $n+1$-gon $\phi_{n+1}^0$. One finds
$$\frac{\phi_{n+1}}{\phi_{n+1}^0} = \frac{ \left((n-1) \cos \frac{\pi}{n} +1 \right) \left( n+1 \right) \tan \left( {\frac {\pi }{n+1}} \right)  \left( 1-\frac{2}{ \left( n+1 \right) } \right) ^{\frac{n+1}{2}}
}{\sin  {\frac {\pi }{n}}   \left(  \left( n-1 \right) +\frac{1}{\cos  {\frac {\pi }{2n}} } 
 \right) ^{2} \left( 1-\frac{2 }{  \left( n-1 \right) +\frac{1} {\left( \cos
 {\frac {\pi }{2n}} \right)    }}
 \right) ^\frac{n-1}{2} \left( 1-\frac{2\,\sin {\frac {\pi }{2n}}}{
   \left( \sin  {\frac {\pi }{n}}   \right) 
 \left(  \left( n-1 \right) +\frac{1}{ \left( \cos {\frac {\pi }{
2n}}   \right)}  \right) } \right)}.$$
Still by {\it Maple} this quotient has only a maximum for $n=6$ and for $n>4$ 
$$\frac{\phi_{n+1}}{\phi_{n+1}^0} > 1.$$
Thus, this example verifies Conjecture (L).\\

{\bf 2.}\ Let us consider now a more general case, which can be viewed as an extension of \cite{ch2}, Example 2 (where $\alpha$ is only closed to $\frac {2\pi }{n}$).  Let a regular $n$-gon with sides of length $a$. We will replace one side by two sides of length respectively $b$ and $c$. We then obtain a $(n+1)$-gon whose $n-1$ sides are of length $a$ one side of length $b$ and one of length $c$. Notice that these two sides are subtended respectively by $\alpha$ and $\frac {\pi }{n}-\alpha, \ 0<\alpha<\frac {2\pi }{n},\ \alpha = \frac {2\pi }{n}$ corresponding to the preceding case where $b=c$. We then compute $b$ and $c$ :
$$a = b\cos \left( \alpha \right) +c\cos \left( -{\frac {\pi }{n}}+\alpha
 \right)\qquad b\sin \left( \alpha \right) =-c\sin \left( -{\frac {\pi }{n}}+\alpha
 \right),$$
we then deduce 
$$ b= a \frac{\sin \left( {\frac {\pi }{n}}-\alpha \right)}{  \left( \sin  {
\frac {\pi }{n}}   \right) }, \qquad c = a \frac{\sin \alpha }{  \left( \sin  {
\frac {\pi }{n}}   \right) }.$$
We then derive the  perimeter, area and pseudo-area of that polygon
$$L_{n+1} = (n-1) a + a \frac{\sin \left( {\frac {\pi }{n}}-\alpha \right)}{  \left( \sin  {
\frac {\pi }{n}}   \right) } + a \frac{\sin \alpha }{  \left( \sin  {
\frac {\pi }{n}}   \right) },$$
$$A_{n+1} = \frac{n}{4}{a}^{2}\cot  {\frac {\pi }{n}}  -\frac{{a}^{2}\sin
 \left( -{\frac {\pi }{n}}+\alpha \right) \sin  \alpha}{2  
 \left( \sin  {\frac {\pi }{n}}   \right)},$$
$$P_{n+1} = \frac{{L_{{n+1}}}^{2}}{4}{ \left( 1-{\frac {2a}{L_{{n+1}}}} \right) ^{\frac{n-1}{2}}\sqrt {1+\frac{2\,a\sin \left( -{\frac {\pi }{n}}+\alpha \right) }{
 \left( \sin  {\frac {\pi }{n}}   \right){L_{{n+1}}
}}}\sqrt {1-\frac{2\,a\sin  \alpha }{  \left( \sin  {                     
\frac {\pi }{n}}   \right) {L_{{n+1}}}}}}.$$
The quotient $\frac{A_{n+1}}{P_{n+1}}$ is 
$$\phi_{n+1} = \frac{ \left( n{a}^{2}\cos  {\frac {\pi }{n}}  +2{a}^{2}\sin \left( {\frac {\pi }{n}}-\alpha \right) \sin  \alpha
   \right)}{  \left( {L_{{n+1}}}^{\frac{-n+3}{2}} \right) 
   \left( L_{{n+1}}-2\,a \right)^{\frac{n-1}{2}}}\times $$ $$\frac{1} {
{\sqrt {(\sin  {\frac {\pi }{n}})\  L_{{n+1}}+2\,a
\sin \left( -{\frac {\pi }{n}}+\alpha \right) }}\ {\sqrt {
(\sin  {\frac {\pi }{n}})\  L_{{n+1}}-2\,a\sin  \alpha
  }}}.$$
	It yields
	$$\phi_{{3}}=\frac{2\,\cos  \alpha  \sin  \alpha}{
   \left( 1+\sin \left( \alpha \right) +\cos \left( \alpha
 \right)  \right) ^{2}}\times $$ $${\frac {1}{\sqrt {1-{\frac {2}{1+\sin
 \left( \alpha \right) +\cos \left( \alpha \right) }}}}}{\frac {1}{
\sqrt {1-{\frac {2\cos \left( \alpha \right) }{1+\sin \left( 
\alpha \right) +\cos \left( \alpha \right) }}}}}{\frac {1}{\sqrt {1-2
\,{\frac {\sin \left( \alpha \right) }{1+\sin \left( \alpha \right) 
+\cos \left( \alpha \right) }}}}} =1,$$
$$\phi_{{4}}={\frac { \left( 3+4\,\cos \left( \frac{\pi}{6} +\alpha
 \right) \sin  \alpha   \right) \sqrt {3}}{4\sqrt {
 \left( \sqrt {3}+\sin  \alpha  -\cos \left( \frac{\pi}{6} +
\alpha \right)  \right) }\sqrt { \left( \sqrt {3}-\sin  \alpha
  +\cos \left( \frac{\pi}{6} +\alpha \right)  \right) } \left( \sin
  \alpha  +\cos \left( \frac{\pi}{6} +\alpha \right)  \right) 
}} =1.$$
Moreover, since \ $0<\alpha<\frac {2\pi }{n}$\ then  
$$lim_{n \rightarrow\infty} \phi_{n+1} =\frac{-\pi \,{{\rm e}}}{f(\alpha)}$$
where
$$f(\alpha)=  -\sqrt {{\frac {\pi +2\,\sin  \alpha
  }{\pi }}}\sqrt {-{\frac {-\pi +2\,\sin  \alpha  
}{\pi }}}{\pi }^{2}+$$ $$2\,\sqrt {{\frac {\pi +2\,\sin  \alpha
  }{\pi }}}\sqrt {-{\frac {-\pi +2\,\sin  \alpha  
}{\pi }}} \left( \cos  \alpha   \right) ^{2}-2\,\sqrt {{
\frac {\pi +2\,\sin  \alpha  }{\pi }}}\sqrt {-{\frac {-
\pi +2\,\sin  \alpha  }{\pi }}}+$$ $$2\,\sqrt {{\frac {\pi +2
\,\sin  \alpha  }{\pi }}}\sqrt {-{\frac {-\pi +2\,\sin
  \alpha  }{\pi }}} \left( \sin  \alpha  
 \right) ^{2}   =-\sqrt {{\pi }^{2}-4\, \left( \sin  \alpha   \right) ^{2}
}\pi.$$
Then

$$lim_{n \rightarrow\infty} \phi_{n+1} = lim_{n \rightarrow\infty} {\frac {{{\rm e}}}{\sqrt {{\pi }^{2}-4\, \left( \sin  \alpha
   \right) ^{2}}}} = \frac{\rm e}{\pi}.$$

By {\it Maple} we may verify for fixed \ $\alpha, \ 0<\alpha<\frac {2\pi }{n}$\ that $\phi_{n+1}$ is strictly decreasing with respect to $n$ from $1$ to $$lim_{n\rightarrow\infty} \phi_{n+1} = \frac{e}{\pi}.$$

{\bf Remark}\ in \cite{ch2}, Proposition 2 one proved for \ $\alpha -\frac{2 \pi}{n} = \epsilon$\ very small
$$\varphi_{n}= \frac{1-\frac{2\varepsilon
^{2}}{(n-2)^{2}}[\frac{(n-2)^{2}}{2n}+\frac{n-2}{2}+\frac{2}{n}-\frac{1}{%
\sin ^{2}\frac{\pi }{n}}]}{{{n \tan{\frac{\pi}{{n}}}}\ (1-{\frac{2}{{n}}}%
)^{n/2}}}.$$
We may prove \ $\varphi_{n}$\ is decreasing and \ $\varphi_{3} = \varphi_{4} =1.$\ Therefore $\varphi_{n}$\ verifies Conjecture (L) : \ $ \phi_n^0 < \phi_n < 1.$\ We hope it is also verified for any \ $0<\alpha<\frac {2\pi }{n}.$\\

\subsection{Generalized Macnab polygons}
{\bf 1.}\  Macnab \cite{m} considered a special polygon with $2n $ sides,  which is a cyclic equiangular
alternate-sided 2n-gon with $n$ sides of length $a$ and $n$ sides of length $b$. The area has been computed by \cite{m} :
$$A_{n,n}=\frac{n}{4\sin \frac{\pi }{n}}[(a^{2}+b^{2})\cos \frac{\pi }{n}
+2ab].$$ 
On the other hand, it is easy to see that the quotient $\varphi _{n,n}= \frac{A_{n,n}}{P_{n,n}}$ may be written
$$\varphi _{n,n}=\frac{[\cos \frac{\pi }{n}+\frac{2ab}{(a+b)^{2}}(1-\cos 
\frac{\pi }{n}]}{n\sin \frac{\pi }{n}[1-\frac{2}{n}+\frac{4ab}{%
n^{2}(a+b)^{2}}]^{\frac{n}{2}}}.$$
Furthermore, we may verify when $n$ growths \ $\varphi _{n,n}$\ is decreasing from $\varphi _{2,2} = 1$ to $\frac{e}{\pi}.$.\\
Therewith in \cite{ch2}, Proposition 1, we proved 
$$\frac{\varphi _{n,n}}{\varphi ^{0}}=\frac{(1-\frac{1}{n})^{n}[\cos \frac{%
\pi }{n}+\frac{2ab}{(a+b)^{2}}(1-\cos \frac{\pi }{n}]}{ 1-\frac{2}{n}+%
\frac{4ab}{n^{2}(a+b)^{2}}]^{\frac{n}{2}}(1+\cos \frac{\pi }{n})} > 1.$$
Thus this polygon verifies Conjecture (L).\\

{\bf 2.}\  Let us consider now another polygon with $k+m$ sides, a polygon alternatively with $k$ sides of length $a$ and $m$ sides of length $b$, $n=k+m$. Such polygon has been considered in \cite{kk}, Theorem 2.3 where its area is computed:
$$L_n = k a+ m b, \quad A_n ={\frac {nab+ \left( k{a}^{2}+m{b}^{2} \right) \cos  
\frac{2\pi}{n}  }{\sin  \frac{2\pi}{n}  }},$$
its pseudo-area is 
$$P_{{n}}=\frac{{L_{{n}}}^{2}}{4} \left( 1-4\,{\frac {a}{L_{{n}}}} \right) ^{\frac{k}{2}
} \left( 1-4\,{\frac {b}{L_{{n}}}} \right) ^{\frac{m}{2}} =$$
$$\frac{\left( 2\,ka+2\,mb \right) ^{2}}{4} \left( 1-{\frac {4 a}{2
\,ka+2\,mb}} \right) ^{\frac{k}{2}} \left( 1-{\frac {4 b}{2\,ka+2\,mb}}
 \right) ^{\frac{m}{2}},$$ and the quotient
$$\phi_n= \frac{ \left( nab+ \left( k{a}^{2}+m{b}^{2} \right) \cos {2
\frac {\pi }{n}}   \right)}{  \left( \sin {\frac {2\pi }
{n}}   \right)  \left( ka+mb \right) ^{2} \left( 
 \left( 1-{\frac {4 a}{2\,ka+2\,mb}} \right) ^{\frac{k}{2}} \right) 
 \left(  \left( 1-{\frac {4 b}{2\,ka+2\,mb}} \right) ^{\frac{m}{2}}
 \right) }.$$

Let us compute now the limit. Notice (6)
$$L_n^2 - 4n \tan \frac{\pi}{n}\ A_n \geq 2 R \tan \frac{\pi}{n}\left[ L_n - 2n R \sin \frac{\pi}{n} \right]$$
implies a bound for the area (here \ $L_n = k a + m b, \ R = \frac{a}{\sin \frac{\pi}{n-k}}$)
$$A_{{n}}<\frac{{L_{{n}}}^{2}}{{4n} \left( \tan  {\frac {\pi }{n}
}   \right) }-\frac{a \left( \frac{L_{{n}}\ \sin \left( {\frac 
{\pi }{n-k}} \right)}{ 2{a}}-n \sin  {\frac {\pi }{n}}  
 \right)}{ {n} \left( \sin \left( {\frac {\pi }{n-k}} \right) 
 \right) },$$ we then derive

$$\phi_n < \frac{\frac{\left( k a+ \left( n-k \right) b \right) ^{2}}{{4n} \left( \tan  {\frac {\pi }{n}
}   \right) }-\frac{a \left( \frac{\left( k a+ \left( n-k \right) b \right) \ \sin \left( {\frac 
{\pi }{n-k}} \right)}{ 2{a}}-n \sin  {\frac {\pi }{n}}  
 \right)}{ {n} \left( \sin \left( {\frac {\pi }{n-k}} \right) 
 \right) }}{\left( k a+ \left( n-k \right) b \right) ^{2} \left( 1-{\frac {2a}{k
a+ \left( n-k \right) b}} \right) ^{\frac{k}{2}} \left( 1-{\frac {2b}{ka+
 \left( n-k \right) b}} \right) ^{\frac{n-k}{2}}}
.$$

 When\ $n \rightarrow \infty$\ the last expression tends to \ $\frac{e}{\pi}.$\\
We then deduce 
$$lim_{n\rightarrow\infty} \phi_{n} \leq \frac{e}{\pi}.$$
By {\it Maple} we may prove \ $\phi_{n}$\ is decreasing when \ $ n \geq 3$\ growths. \\

\subsection{Polygons with bounded sides}\quad However, it is worth considering polygons that are fairly close to the regular one; this will allow for a better test of the conjecture. To do this, let us consider a $n$-gon with perimeter $L_n,$ enclosing an area $A_n$, (\ $P_n$\ is its pseudo-area ) whose sides $a_1, a_2,....,a_n$, satisfying the condition
$$a_i < \frac{k\ L_n}{n}, \ 1<i<n, \ 1 < k < \frac{3}{2}.  \qquad (C_1)$$ 
Since $L_n = \sum_{1<i<n} a_i $ then necessary  $k >1$, \ $k=1$ corresponds to the regular polygon. We prove the following which improves slightly \cite{ch1}, Prop. 4-2.\\

{\bf Proposition 2.1}\ {\it Under the condition $(C_1)$ the quotient \ $\varphi_n = \frac{{A_n}}{{P_n}} $ verifies the following inequalities  
\begin{equation}\frac{4 A_n}{L_n^2 \left( 1-\frac{2}{n}\right)^{\frac{n}{2}}}
 < \varphi_n < \frac{4 A_n}{L_n^2 \left( 1-\frac{2 k}{n}\right)^{\frac{n}{2}}} <\frac{1}{n \tan \frac{\pi}{n} \left(1 - \frac{2 k}{n} \right)^{\frac{n}{2}}},\end{equation}

with equalities if and only if the $n$-gon is regular (corresponding to $k=1$).}\\

{\bf Proof}\quad For the left inequality one proved in \cite{ch1}, Lemma 4-6
$$(1 - {\frac{2a_1}{{L_n}}})(1 - {\frac{2a_2%
}{{L_n}}})(1 - {\frac{2a_3}{{L_n}}}).....(1 - {\frac{2a_n}{{L_n}}}) < \left( 1-\frac{2}{n}\right)^{n},$$
we then deduce 
$$P_n < \frac{L_n^2}{4} \left( 1-\frac{2}{n}\right)^{\frac{n}{2}},$$
$$\frac{4 A_n}{L_n^2 \left( 1-\frac{2}{n}\right)^{\frac{n}{2}}} < \frac{A_n}{P_n}.$$

 Let us prove the right inequality.\\
Indeed, $a_i < \frac{k\ L_n}{n}$ implies $$P_n = {\frac{ L_n^2}{{4}}} \prod_{1<i<n}{\left(1 - {\frac{2a_i}{{L_n}}}\right)^{\frac{1}{2}}} > {\frac{ L_n^2}{{4}}}  \left( 1-\frac{2k}{n} \right)^{\frac{n}{2}}.$$

Therefore $$\phi_n < \frac{A_n}{{\frac{ L_n^2}{{4}}}  \left( 1-\frac{2 k}{n} \right) ^{\frac{n}{2}}}.$$
Moreover, since \ $$A_n < \frac{L^2_n}{4n \tan {\frac {\pi }{n}}}$$ then 
$$\phi_n < \frac{\frac{L^2_n}{4n \tan {\frac {\pi }{n}}}}{{\frac{ L_n^2}{{4}}}  \left( 1-\frac{2 k}{n} \right) ^{\frac{n}{2}}} =
\frac{1}{n \tan \frac{\pi}{n} \left(1 - \frac{2 k}{n} \right)^{\frac{n}{2}}}.$$
On the other hand, for $k=1$, the left and the right sides equal the quotient of the regular $n$-gon. 
$${\varphi_n}^0 = \frac{{A_n}^0}{{P_n}^0} = {\frac{1}{{{n \tan{\frac{\pi}{{n}}}}\ (1-{\frac{2}{{n}}}%
)^{n/2}}}}.$$
We have to prove now the bound for $k$.\ Recall by (6) the necessary condition for the sides \ $a_i, \ 1<i<n$  to be the edge lengths of a cyclic polygon
$$a_i < \frac{L_n}{2}.$$ When the real $k$ is such that \ $ 1 < k < \frac{3}{2},$\ this condition is then satisfied. Indeed, for \ $n \geq 3$ \ the following inequalities hold $$ a_i < \frac{k L_n}{n} < \frac{ L_n}{2}.$$\\

{\bf Remarks 2.2}\\ {\bf 1.}\ This real number $k$ can be used to measure the deviation of a polygon from a regular one with the same area and the same perimeter. The closer $k$ is to $1$, the more nearly regular the polygon becomes. In \cite{ch1}, Prop. 4-2 one proved the inequality
$$\frac{\sqrt {1-{k}^{2} \left( \sin  {\frac {\pi }{n}}  
 \right) ^{2}}}{{n} \left( \sin  {\frac {\pi }{n}}  
 \right)    \left( 1-\frac{2}{n} \right) ^{\frac{n}{2}} } < \varphi_n.$$

For $n$ very large the left side is getting closer to $\frac{e}{\pi}$ :
$$lim_{n\rightarrow \infty}\frac{\sqrt {1-{k}^{2} \left( \sin  {\frac {\pi }{n}}  
 \right) ^{2}}}{{n} \left( \sin  {\frac {\pi }{n}}  
 \right)    \left( 1-\frac{2}{n} \right) ^{\frac{n}{2}} } = \frac{e}{\pi}.$$
We may also notice the equivalence $$\varphi_n = \frac{{A_n}}{{P_n}} = 1 \Leftrightarrow  
 k = \frac{n}{2} - \frac{n}{2 \left( \tan{\frac{\pi}{{n}}} \right)^{\frac{2}{n}}}.$$ In addition, it is easy to see (by {\it Maple}) that $n=3$ or $n=4 \Leftrightarrow k=1$. We then deduce the inequalities since $k$ is increasing with $n$
$$1 < k < \frac{n}{2} - \frac{n}{2 \left( \sin{\frac{\pi}{{n}}} \right)^\frac{2}{n}} < \log \pi < \frac{3}{2}.$$
Notice also for $n=3$ or $n=4$ we have necessarily $k=1$\ (in these cases \ $\phi_3=\phi_4 =1$)\ .\ Indeed, the solution of 
$$\frac{1}{{n} \left( \tan  {\frac {\pi }{n}}   \right) 
 \   \left( 1-2\,{\frac {k}{n}} \right) ^{\frac{n}{2}}  }=1,$$
is $$k = \frac{n}{2}-\frac{n}{2}{\rm e}^{\frac{2\,\ln \frac{1}{ \left( {n} \left( \tan \left( {\frac {
\pi }{n}} \right)  \right)  \right)}}{ {n}}}.$$
For $n=3$ or $n=4$ it yields $k=1$. Otherwise \ $1 < k < \log \pi.$\\

{\bf 2.}\ On the other hand, we may prove by {\it Maple} for $k >1$ and $n> 4$ : $$\frac{1}{\left( 1-\frac{2}{n} \right) ^{\frac{n}{2}} } < \exp\left(\frac{k}{1-\frac{2k}{n}}\right),$$
which implies $$\frac{1}{n \tan \frac{\pi}{n} \left(1 - \frac{2 k}{n} \right)^{\frac{n}{2}}} < \frac{\exp\left(\frac{k}{1-\frac{2k}{n}}\right)}{n \tan \frac{\pi}{n}}.$$
This means Proposition 2.1 improves \cite{ch1}, Prop. 4-2.\\

Moreover, in using various approaches of Bonnesen inequalities it is possible to improve Proposition 2.1. Indeed, we may use the above inequality (6) 
$$L_n^2 = 4  n \tan \frac{\pi}{n} A_n \geq \left(2 \pi n \sin \frac{\pi}{n} -L_n \right)^2,\qquad (6)$$
with equality if and only if the polygon is regular. We then obtain\\ 

{\bf Proposition 2.3}\ {\it Under the condition $(C_1)$ the quotient \ $\varphi_n = \frac{{A_n}}{{P_n}} $ verifies the following inequalities 
$$ \varphi_n < \frac{4 A_n}{L_n^2 \left( 1-\frac{2 k}{n}\right)^{\frac{n}{2}}} \leq  \frac{ 1 - \left(\frac{n \tan \frac{\pi}{n} A_n}{L_n}-1 \right)^2}{n\tan \frac{\pi}{n} \left( 1-\frac{2 k}{n}\right)^{\frac{n}{2}}}
 <\frac{1}{n \tan \frac{\pi}{n} \left(1 - \frac{2 k}{n} \right)^{\frac{n}{2}}},$$ 
 with equalities if and only if the $n$-gon is regular (corresponding to $k=1$).}\\

{\bf Proof }\ Indeed, (6) is equivalent to 
$$\frac{4  n \tan \frac{\pi}{n} A_n}{L_n^2} \leq 1 - \left(\frac{n \tan \frac{\pi}{n} A_n}{L_n}-1 \right)^2,$$
which implies 
$$\phi_n = \frac{A_n}{P_n} < \frac{4 A_n}{L_n^2 \left( 1-\frac{2 k}{n}\right)^{\frac{n}{2}}} \leq  \frac{1 - \left(\frac{n \tan \frac{\pi}{n} A_n}{L_n}-1 \right)^2}{n\tan \frac{\pi}{n} \left( 1-\frac{2 k}{n}\right)^{\frac{n}{2}}} < \frac{1}{n \tan \frac{\pi}{n} \left(1 - \frac{2 k}{n} \right)^{\frac{n}{2}}}.$$  

\section{Concluding remark}

To conclude this paper, let us notice that we can provide other examples of polygons, or
even arbitrary polygon as models to obtain various forms of functions $\phi_n$, in order to verify Conjecture $(L)$ and to
establish some other types of analytic and geometric isoperimetric inequalities. Of course, intensive use of {\it Maple} or any other tool is required.\\
The idea would be to expand on the final section regarding bounds on the lengths. It would be wise to highlight other, more general bounds than the one given by $(C_1)$ : 
$$a_i < \frac{k\ L_n}{n}, \ 1<i<n, $$ that ensure the polygon remains close to a regular one; these would, of course, depend on one or more parameters.
Nevertheless, as is evident, fully resolving this conjecture remains a very difficult task. However, partial answers can be provided.\\


\begin{thebibliography}{12}

\bibitem{b} M. Berger\ {\it Polygons, polyhedra, polytopes},  Geometry Revealed. Springer, Berlin, Heidelberg, p.505-561, (2010). 10.1007/978-3-540-70997-8-8.

\bibitem{ch1} R. Chouikha\ {\it Probleme de P. Levy sur les polygones articules}
\newline
C. R. Math. Report, Acad of Sc. of Canada, vol 10, p. 175-180, 1988. 

\bibitem{ch2} R. Chouikha\ {\it Problems on polygons and Bonnesen-type inequalities}, Indag., Vol. 10, Issue 4, 1999, p. 495-506.\\
10.1016/S0019-3577(00)87902-1

\bibitem{dl} P. Dulio, E. Laeng\ {\it Generalization of Heron’s and Brahmagupta’s equalities to any
cyclic polygon}\qquad Aequat. Math. 95 (2021), 941–952. \\ doi.org/10.1007/s00010-020-00771-w.

\bibitem{kk} H.T. Ku, M.C. Ku\ {\it Analytic Isoperimetric Inequalities}, Math. Ineq. Appl., Volume 3, Number 4 (2000), 459–472.

\bibitem{kkz} H.T. Ku, M.C. Ku, X.M. Zhang\ {\it Analytic and geometric
isoperimetric inequalities} \qquad J. of Geometry, vol 53, p.100-121, 1995.


\bibitem{l} P. Levy\ {\it Le probleme des isoperimetres et des polygones articules}
\newline
Bull. Sc. Math., 2eme serie, 90, p.103-112, 1966.

\bibitem{m}  D. S. Macnab {\it Cyclic polygons and related questions}, Math. Gazette 65 (1981), 22–28.

\bibitem{o1} R. Osserman\ {\it The isoperimetric inequalities}\newline
Bull. Amer. Math. Soc., vol 84, p.1182-1238, 1978.

\bibitem{o2} R. Osserman\ {\it Bonnesen-style isoperimetric inequalities}\newline
Amer. Math. Monthly, vol 1, p. 1-29, 1979.

\bibitem{pak}  I. Pak,\ {\it  The area of cyclic polygons: Recent progress on Robbins’ conjectures},\quad Adv.
in Appl. Math., Vol. 34, Issue 4, 2005, p. 690–696.

 \bibitem{pech} P. Pech, {\it Computations of the Area and Radius of Cyclic Polygons Given by the Lengths
of Sides},\quad ADG2004 (Hong, H. and Wang, D., eds.), LNAI, 3763, Gainesville, Springer,
2006, 44–58.

\bibitem{p} I. Pinelis, \ {\it Cyclic polygons with given edge lengths}, \ J. Geom. 00 (2005), P.1-16. 
DOI 10.1007/s00022-005-1752-8.

\bibitem{rob}  D.P. Robbins, {\it Areas of polygons inscribed in a circle}\quad Discrete Comput. Geom. 12(2),
223–236 (1994).

\bibitem{sv} D. Svrtan\ {\it On circumradius equations of cyclic polygons}\\ Proc. of the 4th Croatian Combin. Days
Sept. 22 – 23, 2022

\bibitem{zd} C. Zeng, X. Dong\ {\it On Some Discrete Bonnesen-style Isoperimetric Inequalities}, \  Acta Math. Sinica, Vol 41, p. 1447–1461, (2025) doi.org/10.1007/s10114-025-3281-8.

\bibitem{z} X.M. Zhang\ {\it Bonnesen-style inequalities and pseudo-perimeters for
polygons}\qquad J. of Geometry, vol 60, p.188-201, 1997.\newline


\end{thebibliography}
\end{document}